\documentclass[11pt]{article}
\usepackage[margin=1.3in]{geometry}
\usepackage[T1]{fontenc}
\usepackage{lmodern}
\usepackage{amsmath,amssymb,amsthm,mathtools}
\usepackage{enumitem}
\usepackage[colorlinks=true,linkcolor=blue,citecolor=blue,urlcolor=blue]{hyperref}
\usepackage[nameinlink,noabbrev,capitalise]{cleveref}
\hypersetup{pdftitle={Isolated singularities of weak solutions to the quasilinear Lane-Emden equation},
  pdfauthor={Meiqing Xu and Hui Yang}}
\newtheorem{theorem}{Theorem}[section]
\newtheorem{proposition}[theorem]{Proposition}
\newtheorem{lemma}[theorem]{Lemma}

\theoremstyle{remark}

\newcommand{\R}{\mathbb{R}}

\newcommand{\Sph}{\mathbb{S}}
\newcommand{\dd}{\,\mathrm d}
\newcommand{\loc}{\mathrm{loc}}
\newcommand{\K}{K_{n,p,q}}
\DeclareMathOperator{\diver}{div}

\setlist[enumerate]{label=\textup{(\roman*)},leftmargin=2em}
\title{\bfseries Isolated singularities of weak solutions to the quasilinear Lane-Emden equation}
\author{Meiqing Xu \and Hui Yang}
\date{}

\begin{document}
\maketitle

\begin{abstract}
We study isolated singularities of nonnegative weak solutions to
\[
-\Delta_pu=u^q\quad\text{in }B_1\setminus\{0\},
\]
where $\Delta_p$ is the $p$-Laplacian, $1<p<n$, and
$n(p-1)/(n-p)<q<np/(n-p)-1$.
We prove that either the origin is a removable singularity or
$u(x)\sim K_{n,p,q}|x|^{-p/(q+1-p)}$ as $x\to0$, where
$K_{n,p,q}>0$ is an explicit constant depending only on $n,p,q$.
The result holds without radial symmetry or an a priori growth assumption.
The proof combines a removability criterion with a rigidity theorem for
positive solutions in $\mathbb R^n\setminus\{0\}$ satisfying two-sided
power bounds. To handle the infinite $p$-Dirichlet energy at the singularity,
we construct a continuous Steiner deformation that fixes all sufficiently
high superlevel sets. The resulting energy comparison yields radial symmetry
of the blow-up limits.


\end{abstract}

\section{Introduction}
In this paper, we study the asymptotic behavior near the origin of nonnegative weak
solutions to
\begin{equation}\label{eq:main}
 -\Delta_pu=u^q\quad\text{in }B_1\setminus\{0\},
\end{equation}
where $\Delta_pu=\diver\bigl(|\nabla u|^{p-2}\nabla u\bigr)$ is the classical $p$-Laplacian, $1<p<n$, $B_1\subset\R^n$ is the unit ball centered at the origin and
\begin{equation}\label{eq:range}
 \frac{n(p-1)}{n-p}<q<\frac{np}{n-p}-1.
\end{equation}
The two endpoints in \eqref{eq:range} are the Serrin and Sobolev critical
exponents. In this range, equation \eqref{eq:main} admits the singular
solution
\begin{equation}\label{eq:singular-profile}
 \Phi(x)=\K|x|^{-\alpha},\qquad
 \alpha=\frac{p}{q+1-p},\qquad
 \K=\bigl[\alpha^{p-1}\bigl(n-p-(p-1)\alpha\bigr)\bigr]^{\frac1{q+1-p}}.
\end{equation}
Our purpose is to show that $\Phi$ describes the leading behavior of
every solution with a nonremovable isolated singularity.

For $p=2$ and $n\ge3$, isolated singularities of $-\Delta u=u^q$ have a classical
classification. Below the Serrin exponent, $1<q<n/(n-2)$, the nonremovable
alternative has the fundamental-solution form $u(x)\sim a|x|^{2-n}$ with
$a>0$ depending on the solution; see Br\'ezis--Lions~\cite{BrezisLions1981} and Lions~\cite{Lions1980}. In the range $n/(n-2)<q<(n+2)/(n-2)$, Gidas--Spruck~\cite{GidasSpruck1981} established the precise power
asymptotics with a coefficient determined by the equation. At the Serrin exponent, Aviles~\cite{Aviles1983,Aviles1987} obtained a
logarithmic correction to the fundamental solution for a nonremovable singularity. At the Sobolev
exponent, Caffarelli--Gidas--Spruck~\cite{CaffarelliGidasSpruck1989}
proved asymptotic radial symmetry and described the possible radial
profiles. These include Fowler solutions, whose normalized profiles
are periodic in the logarithmic radial variable. Further results on
asymptotic symmetry and refined asymptotics can be found
in~\cite{KorevaarMazzeoPacardSchoen1999,Li1996}.

For quasilinear equations, the study of isolated
singularities goes back to Serrin~\cite{Serrin1964,Serrin1965}.
Precise asymptotics for singular $p$-harmonic functions were obtained
by Kichenassamy--V\'eron~\cite{KichenassamyVeron1986}.
For the equation \eqref{eq:main}, Guedda--V\'eron~\cite{GueddaVeron1988} classified the isolated singularities of positive radial
solutions in the subcritical Sobolev range. Bidaut-V\'eron~\cite{BidautVeron1989} obtained the fundamental-solution asymptotics for \eqref{eq:main} below the Serrin threshold without assuming radial symmetry. The two critical endpoints have also attracted recent attention. Ma--Zang~\cite{MaZang2025} studied the isolated singularities of a Serrin-critical quasilinear inequality. Xu and Yang~\cite{XuYang2026} established exact logarithmic asymptotics for nonnegative weak solutions of \eqref{eq:main} in the Serrin critical case; see also
Ma and Zang~\cite{MaZang2026} for corresponding results in the class of positive $C^2$ solutions. At the Sobolev critical exponent, Yang~\cite{YangCritical2026} proved that every solution of \eqref{eq:main} with a nonremovable singularity is asymptotic to
a single radial Fowler solution. In this paper, we study the intermediate range between the Serrin and Sobolev critical exponents.

A nonnegative \emph{weak} solution of \eqref{eq:main} is a function
$u\ge0$ with $u\in W^{1,p}_{\loc}(B_1\setminus\{0\})$ such that
\begin{equation}\label{eq:weak-formulation}
 \int_{B_1\setminus\{0\}}
 |\nabla u|^{p-2}\nabla u\cdot\nabla\varphi\dd x
 =
 \int_{B_1\setminus\{0\}}u^q\varphi\dd x
\end{equation}
for every $\varphi\in C_c^\infty(B_1\setminus\{0\})$.
In the subcritical range the right-hand side is locally integrable,
and such a solution has a locally $C^{1,\eta}$ representative in $B_1\setminus\{0\}$ for
some $\eta>0$. We recall this regularity argument in
Section~\ref{sec:prelim}. The singularity at the origin is called \emph{removable} if $u$ extends
to a locally bounded weak solution of \eqref{eq:main} throughout $B_1$.
By standard local regularity, this extension belongs to $C^{1,\eta}_{\loc}(B_1)$ for some $\eta\in(0,1)$.

By the universal estimate of
Serrin--Zou~\cite[Theorem~IV(b)]{SerrinZou2002},
every nonnegative solution of \eqref{eq:main} under \eqref{eq:range} satisfies
\begin{equation}\label{eq:upper}
 u(x)\le C_0|x|^{-\frac{p}{q+1-p}}, \qquad 0<|x|<1/2,
\end{equation}
where $C_0=C_0(n,p,q)>0$ is independent of $u$. It remains to characterize removable
singularities and determine the asymptotic behavior in the nonremovable
case, without assuming radial symmetry. Our main result is the following.

\begin{theorem}\label{thm:classification}
Assume that \eqref{eq:range} holds. Let $u$ be a nonnegative weak solution of \eqref{eq:main}. Then exactly one of the following alternatives holds:
\begin{enumerate}
 \item The singularity at the origin is removable.
 \item As $r\downarrow0$,
 \begin{equation}\label{eq:C1limit}
 r^\alpha u(r\,\cdot) \to \K|\cdot|^{-\alpha}
 \quad\text{in }C^1_{\loc}(\R^n\setminus\{0\}),
 \end{equation}
 where $\alpha$ and $\K$ are given by \eqref{eq:singular-profile}. In particular,
\begin{equation}\label{eq:precise}
 \lim_{x\to0}|x|^\alpha u(x)=\K.
 \end{equation}
\end{enumerate}
\end{theorem}

The first difficulty in proving Theorem~\ref{thm:classification} is that, unlike in the
semilinear case \cite{GidasSpruck1981,YangZou}, the derivative of the natural energy on the cylinder no longer has a definite sign. Indeed, for $v(t,\theta)=e^{-\alpha t}u(e^{-t}\theta)$, set
$z=v_t+\alpha v$ and $Q=(z^2+|\nabla_{\Sph^{n-1}}v|^2)^{1/2}$. The cylinder equation is
\begin{equation}\label{eq:cylinder}
 (Q^{p-2}z)_t-\beta Q^{p-2}z
 +\diver_{\Sph^{n-1}}(Q^{p-2}\nabla_{\Sph^{n-1}}v)+v^q=0.
\end{equation}
For smooth profiles, the natural corrected energy is
\[
 \mathcal H(t)=\int_{\Sph^{n-1}}\left[
 Q^{p-2}zv_t-\frac{Q^p}{p}+\frac{v^{q+1}}{q+1}
 +\frac{\delta\alpha^{p-1}}p v^p\right]\dd\theta.
\]
Formal differentiation and the cylinder equation give
\begin{equation}\label{eq:no-sign}
 \mathcal H'(t)=-\delta\int_{\Sph^{n-1}}
 \bigl[Q^{p-2}z-(\alpha v)^{p-1}\bigr]v_t\dd\theta,
\end{equation}
where $\delta=p(\alpha+1)-n >0$. For $p=2$ this is $-\delta\int v_t^2$, and hence $\mathcal H$ is monotone; but for general $p$ the angular
gradient prevents that conclusion. Instead of such a monotonicity formula for $p=2$, we shall use a two-power barrier and Picone's identity to prove that the singularity at the origin is removable if $\liminf\limits_{x\to0}|x|^\alpha u(x)=0$, and thereby classify the isolated singularities. A standard rescaling argument then reduces Theorem~\ref{thm:classification} to the classification of positive
solutions $U$ in $\R^n\setminus\{0\}$ with the same two-sided bounds.
Rigidity for already homogeneous profiles was studied by
Lin--Ma~\cite{LinMa2024}. The punctured-space symmetry result of
Esposito--Muglia--Sciunzi--Vuono~\cite{EspositoMugliaSciunziVuono2026}
assumes $U^{q+1-p}\in L^{n/p}$ outside a ball. In our setting,
$U^{q+1-p}$ is comparable to $|x|^{-p}$, so this integrability condition is unavailable.

The second main difficulty is to establish the radial symmetry of these limiting solutions in $\R^n\setminus\{0\}$. For general  $p\ne2$, the classical method of moving planes cannot be applied directly due to the lack of a Kelvin transformation \cite{CaffarelliGidasSpruck1989}. Motivated by Brock's continuous Steiner symmetrization~\cite{Brock1995,Brock2000} and its application to the critical problem in~\cite{YangCritical2026}, we establish radial symmetry via Steiner symmetrization. However, in the subcritical regime, the stronger energy divergence near the singularity
prevents a direct use of the logarithmic cutoff argument in~\cite{YangCritical2026}, so we need to construct a Steiner flow to keep all
sufficiently high superlevel sets fixed. The resulting variation
vanishes near the singularity and can be used in the weak equation.
Equimeasurability makes its energy loss of order $O(t^2)$, whereas
a failure of local symmetry produces a loss of order $t$. This gives
local symmetry in every direction. The positivity of the source
then yields radial symmetry, and a radial energy identity identifies
the unique profile \eqref{eq:singular-profile}.

The paper is organized as follows. Section~\ref{sec:prelim} collects
the regularity, scaling and extension properties used in the proof.
In Section~\ref{sec:removability}, we establish the removability
criterion. Section~\ref{sec:steiner} contains the construction and
energy estimates for the modified Steiner deformation.
In Section~\ref{sec:rigidity}, we prove the punctured-space
classification and complete the proof of
Theorem~\ref{thm:classification}.
Throughout the paper, $C$ denotes a positive constant that may
change from line to line; its dependence is indicated when needed.

\section{Preliminaries}\label{sec:prelim}
We first record the estimates used in the proof. Let $u\ge0$ be a weak
solution of \eqref{eq:main}. On every ball compactly contained in
$B_1\setminus\{0\}$, Sobolev embedding gives $u\in L^{p^*}$, where
$p^*=np/(n-p)$. Thus
\[
 -\Delta_pu=a(x)u^{p-1},\qquad
 a=u^{q+1-p}\in L^t_{\loc},\qquad
 t=\frac{p^*}{q+1-p}>\frac np.
\]
Moreover, $q<p^*-1$ allows the weak equation to be tested with
$W^{1,p}_0$ functions on such balls. The local boundedness estimate
of Serrin~\cite[Theorem~1]{Serrin1964} and the regularity results
in~\cite{DiBenedetto1983,Tolksdorf1984} imply that
$u\in C^{1,\eta}_{\loc}(B_1\setminus\{0\})$ for some $\eta\in(0,1)$.
By the strong maximum principle, either $u\equiv0$ or $u>0$.
In what follows, we consider the latter case whenever positivity is needed.

For later use, set
\begin{equation}\label{eq:parameters}
 \beta=\frac{(n-p)q-n(p-1)}{q+1-p},
 \qquad \delta=\frac{np-(n-p)(q+1)}{q+1-p}.
\end{equation}
The range \eqref{eq:range} gives
\begin{equation}\label{eq:parameter-signs}
 \frac{n-p}{p}<\alpha<\frac{n-p}{p-1},
 \qquad\beta>0,\quad\delta>0.
\end{equation}
For $r>0$, define
\begin{equation}\label{eq:scaling}
 u_r(y)=r^\alpha u(ry).
\end{equation}
Then $u_r$ satisfies the same equation in $B_{1/r}\setminus\{0\}$.
By \eqref{eq:upper} and interior regularity in~\cite{DiBenedetto1983,Tolksdorf1984}, there exists $\eta\in(0,1)$
such that the family $\{u_r\}$ is uniformly bounded in $C^{1,\eta}$ on each
fixed compact annulus, for all sufficiently small $r$.
In particular, it is precompact in $C^1$ on such annuli, and
\begin{equation}\label{eq:gradient-bound}
 |\nabla u(x)|\le C|x|^{-\alpha-1}
 \qquad(0<|x|<1/4).
\end{equation}
Since $u_r^{q+1-p}$ is uniformly bounded on fixed annuli, Harnack's
inequality gives, for every $L>1$,
\begin{equation}\label{eq:annular-harnack}
 \sup_{r\le |x|\le Lr}u(x)
 \le H_L\inf_{r\le |x|\le Lr}u(x),
\end{equation}
where $H_L$ is independent of sufficiently small $r$.

The equation extends across the origin in the distributional sense.
\begin{lemma}\label[lemma]{lem:no-mass}
Let $u\ge0$ be a weak solution of \eqref{eq:main}. Then
$u^q,\,|\nabla u|^{p-1}\in L^1(B_{1/4})$. Moreover,
\[
 \int_{B_{1/4}}|\nabla u|^{p-2}\nabla u\cdot\nabla\varphi\dd x
 =\int_{B_{1/4}}u^q\varphi\dd x
 \qquad\text{for every }\varphi\in C_c^\infty(B_{1/4}).
\]
\end{lemma}

\begin{proof}
The asserted integrability follows from \eqref{eq:upper} and
\eqref{eq:gradient-bound}. Choose a smooth cutoff
$\chi_\rho$ such that $\chi_\rho=0$ in $B_\rho$,
$\chi_\rho=1$ outside $B_{2\rho}$, and
$|\nabla\chi_\rho|\le C/\rho$.
Testing the punctured equation with
$\chi_\rho\varphi$, we have
\[
\begin{aligned}
&\int_{B_{1/4}}\chi_\rho
|\nabla u|^{p-2}\nabla u\cdot\nabla\varphi\dd x+
\int_{B_{2\rho}\setminus B_\rho}
\varphi|\nabla u|^{p-2}\nabla u
\cdot\nabla\chi_\rho\dd x
=\int_{B_{1/4}}u^q\chi_\rho\varphi\dd x.
\end{aligned}
\]
Since $0\le\chi_\rho\le1$ and
$\chi_\rho\to1$ almost everywhere in $B_{1/4}$,
the integrability of $u^q$ and $|\nabla u|^{p-1}$
allows us to apply the dominated convergence theorem:
\[
\begin{aligned}
&\lim_{\rho\to0}
\int_{B_{1/4}}\chi_\rho
|\nabla u|^{p-2}\nabla u\cdot\nabla\varphi\dd x
=
\int_{B_{1/4}}
|\nabla u|^{p-2}\nabla u\cdot\nabla\varphi\dd x,\\
&\lim_{\rho\to0}
\int_{B_{1/4}}u^q\chi_\rho\varphi\dd x
=
\int_{B_{1/4}}u^q\varphi\dd x.
\end{aligned}
\]
The remaining term  can be estimated by
\[
 \left|\int_{B_{2\rho}\setminus B_\rho}
 \varphi|\nabla u|^{p-2}\nabla u\cdot\nabla\chi_\rho\dd x\right|
 \le C\rho^{-1}\int_{B_{2\rho}\setminus B_\rho}
 |\nabla u|^{p-1}\dd x
 \le C\rho^\beta\longrightarrow0.
\]
 Letting $\rho\to0$, we conclude that
\[
\int_{B_{1/4}}
|\nabla u|^{p-2}\nabla u\cdot\nabla\varphi\dd x
=
\int_{B_{1/4}}u^q\varphi\dd x
\]
for every $\varphi\in C_c^\infty(B_{1/4})$.
\end{proof}

\section{Removability of the singularity}\label{sec:removability}
In this section, we prove that the origin is removable whenever
$\liminf_{x\to0}|x|^\alpha u(x)=0$. We begin with a smallness criterion.

\begin{lemma}\label[lemma]{lem:rem-smallness}
There exists a smallness constant
$\varepsilon_*=\varepsilon_*(n,p,q)>0$ such that for any $R>0$, let $u\in W^{1,p}_{\loc}(B_R\setminus\{0\})$ be a nonnegative weak solution of $-\Delta_pu=u^q$ in $B_R\setminus\{0\}$.
If
\begin{equation}\label{eq:rem-small-hypothesis}
 u(x)\le\varepsilon_*|x|^{-\alpha}
 \qquad(0<|x|<R),
\end{equation}
then $u$ extends to a weak solution in $B_R$, and the extension belongs to $C^{1,\vartheta}_{\mathrm{loc}}(B_R)$ for some $\vartheta\in(0,1)$.
\end{lemma}
Although the constant $C_0$ in \eqref{eq:upper} is also
independent of the solution, that estimate alone does not
ensure the smallness condition \eqref{eq:rem-small-hypothesis},
which requires the normalized solution $|x|^\alpha u(x)$
to remain below the prescribed constant $\varepsilon_*$.
\begin{proof}
We may assume that $u>0$ and, by reducing $R$, that $u$ is
$C^1$ near $\partial B_R$. Fix
\begin{equation}\label{eq:rem-auxiliary-powers}
 0<\gamma<\min\left\{\alpha,\frac{n-p}{p},\frac pq\right\},
 \qquad \alpha<\Gamma<\frac{n-p}{p-1},
\end{equation}
and put
\[
 b_\Gamma=n-p-(p-1)\Gamma.
\]
Choose $\varepsilon_*$ such that $\varepsilon_*^{q+1-p}<b_\Gamma\gamma^{p-1}$.

For $A,B\ge0$, not both zero, let
$w(r)=Ar^{-\gamma}+Br^{-\Gamma}$. A direct computation gives
\begin{align}
 -\Delta_pw
 &=r^{-p}(\gamma Ar^{-\gamma}+\Gamma Br^{-\Gamma})^{p-2}
 \times\bigl[\gamma\{n-p-(p-1)\gamma\}Ar^{-\gamma}
                  +\Gamma b_\Gamma Br^{-\Gamma}\bigr]
 \notag\\
 &\ge b_\Gamma r^{-p}(\gamma Ar^{-\gamma}+\Gamma Br^{-\Gamma})^{p-1}
 \ge b_\Gamma\gamma^{p-1}r^{-p}w^{p-1}.
 \label{eq:rem-two-power-supersolution}
\end{align}
For $0<\rho<R$, define
\[
 w_\rho(r)=\varepsilon_*R^{\gamma-\alpha}r^{-\gamma}
          +\varepsilon_*\rho^{\Gamma-\alpha}r^{-\Gamma}.
\]
Then $w_\rho\ge u$ on $\partial(B_R\setminus\overline B_\rho)$.
To prove the same inequality in the annulus, set
$E=\{u>w_\rho\}$ and $\zeta=(u^p-w_\rho^p)_+$. The functions
$\zeta/u^{p-1}$ and $\zeta/w_\rho^{p-1}$ belong to
$W^{1,p}_0(B_R\setminus\overline B_\rho)$.
Testing the equation for $u$ with
$\zeta/u^{p-1}$, we obtain
\[
 \int_E |\nabla u|^{p-2}\nabla u \cdot\nabla\left(\frac{\zeta}{u^{p-1}}\right)\dd x =\int_E u^{q+1-p}\zeta\dd x \le \varepsilon_*^{q+1-p} \int_E |x|^{-p}\zeta\dd x.
\]
On the other hand, testing the supersolution
inequality \eqref{eq:rem-two-power-supersolution}
with $\zeta/w_\rho^{p-1}$ gives
\[
 \int_E |\nabla w_\rho|^{p-2}\nabla w_\rho \cdot\nabla\left(\frac{\zeta}{w_\rho^{p-1}}\right)\dd x \ge b_\Gamma\gamma^{p-1} \int_E |x|^{-p}\zeta\dd x.
\]
Picone's inequality
\cite{AllegrettoHuang1998} states that
\[
 P(v,w):=|\nabla v|^p-
 |\nabla w|^{p-2}\nabla w\cdot\nabla(v^p/w^{p-1})\ge0.
\]
Subtracting the two tested inequalities and
using Picone's inequality, we therefore obtain
\[
 0 \le\int_E \bigl[P(u,w_\rho)+P(w_\rho,u)\bigr]\dd x \le \bigl(\varepsilon_*^{q+1-p} -b_\Gamma\gamma^{p-1}\bigr) \int_E |x|^{-p}(u^p-w_\rho^p)\dd x.
\]
By our choice of $\varepsilon_*$,
\[
\varepsilon_*^{q+1-p}
-b_\Gamma\gamma^{p-1}<0.
\]
Since $u^p-w_\rho^p>0$ on $E$,
the preceding inequality forces $E$ to be empty.
Thus $u\le w_\rho$ in
$B_R\setminus\overline B_\rho$. Letting $\rho\downarrow0$ yields
\begin{equation}\label{eq:rem-improved-growth}
 u(x)\le \varepsilon_*R^{\gamma-\alpha}|x|^{-\gamma}.
\end{equation}
For the rescaling $v_r(y)=r^\gamma u(ry)$, we have
\[
 -\Delta_pv_r=r^{p-\gamma(q+1-p)}v_r^q.
\]
Since $\gamma<\alpha$, the coefficient on the right is bounded
as $r\downarrow0$. Estimate \eqref{eq:rem-improved-growth} and
interior gradient estimates therefore give
\begin{equation}\label{eq:rem-improved-gradient}
 |\nabla u(x)|\le C|x|^{-\gamma-1}.
\end{equation}
Since $p(\gamma+1)<n$, estimates
\eqref{eq:rem-improved-growth} and
\eqref{eq:rem-improved-gradient} imply that
$u$ and $\nabla u$ belong to $L^p(B_R)$,
where $\nabla u$ is initially understood on the
punctured ball. Let $\chi_\rho$ be the cutoff
function used in Lemma~\ref{lem:no-mass}.
Since $u\chi_\rho\in W^{1,p}(B_R)$, we have
\[
\begin{aligned}
\|u\chi_\rho-u\|_{L^p(B_R)}^p
\le \int_{B_{2\rho}}|u|^p\dd x\le C\rho^{n-p\gamma}\longrightarrow0.
\end{aligned}
\]
Moreover, using
$\nabla(u\chi_\rho)
=\chi_\rho\nabla u+u\nabla\chi_\rho$,
we obtain
\[
\begin{aligned}
\|\nabla(u\chi_\rho)-\nabla u\|_{L^p(B_R)}^p\le C\int_{B_{2\rho}}|\nabla u|^p\dd x
+C\int_{B_{2\rho}\setminus B_\rho}
|u\nabla\chi_\rho|^p\dd x\le C\rho^{n-p(\gamma+1)}
\longrightarrow0.
\end{aligned}
\]
Thus $\{u\chi_\rho\}$ is a Cauchy family in
$W^{1,p}(B_R)$. By completeness, it converges
to some $v\in W^{1,p}(B_R)$. Since
$u\chi_\rho\to u$ in $L^p(B_R)$, we have
$v=u$ almost everywhere. Consequently,
$u\in W^{1,p}(B_R)$, proving the Sobolev
extension across the origin.

Arguing as in the proof of Lemma~\ref{lem:no-mass},
we test the equation with $\varphi\chi_\rho$
and let $\rho\downarrow0$.
By \eqref{eq:rem-improved-gradient}, the cutoff
error is bounded by
$C\rho^{n-p-(p-1)\gamma}\to0$.
Thus the weak equation holds in the whole ball $B_R$.

Finally, $\gamma q<p$ allows us to choose
$n/p<t<n/(\gamma q)$, so that $u^q\in L^t(B_R)$.
The local boundedness theorem for $p$-Poisson equations
\cite[Theorem~1]{Serrin1964} gives
$u\in L^\infty_{\loc}(B_R)$. The claimed regularity $C^{1,\vartheta}_{\mathrm{loc}}(B_R)$ now follows
from~\cite{DiBenedetto1983,Tolksdorf1984}.
\end{proof}

\begin{proposition}\label[proposition]{prop:removability}
Let $u\ge0$ be a weak solution of \eqref{eq:main}. The origin is
removable if and only if
\[
 \liminf_{x\to0}|x|^\alpha u(x)=0.
\]
If the origin is nonremovable, there exist $c_0>0$ and
$r_0\in(0,1/2)$ such that
\begin{equation}\label{eq:two-sided-local}
 c_0|x|^{-\alpha}\le u(x)\le C_0|x|^{-\alpha}
 \qquad(0<|x|<r_0),
\end{equation}
where $C_0=C_0(n,p,q)$ is the constant in \eqref{eq:upper}.
\end{proposition}

\begin{proof}
It suffices to consider $u>0$. Set
\[
 M(r)=\max_{|x|=r}|x|^\alpha u(x),\qquad 0<r<1/2.
\]
Suppose that $|x_j|^\alpha u(x_j)\to0$ for some sequence $x_j\to0$.
By \eqref{eq:annular-harnack}, $M(|x_j|)\to0$.
We claim that $M(r)\to0$ as $r\downarrow0$.

Suppose otherwise. Choose
\[
 0<\varepsilon<\min\{\varepsilon_*,\limsup_{r\downarrow0}M(r)\}.
\]
For all sufficiently large $j$, let $R_j>|x_j|$ be the first radius
at which $M$ reaches $\varepsilon$. Such a radius exists by continuity. Then
\begin{equation}\label{eq:rem-crossing}
 M(R_j)=\varepsilon,
 \qquad M(r)\le\varepsilon\quad(|x_j|\le r\le R_j).
\end{equation}
 Moreover, $R_j\to0$: given
$a>0$, choose $r_a\in(0,a)$ with $M(r_a)>\varepsilon$; for large
$j$, $|x_j|<r_a$ and hence $R_j\le r_a$.
We also have
\begin{equation}\label{eq:rem-crossing-separation}
 \frac{R_j}{|x_j|}\longrightarrow\infty.
\end{equation}
Indeed, if $R_j/|x_j|\le L$ along a subsequence, then
\eqref{eq:annular-harnack} gives
\[
 \varepsilon=M(R_j)\le L^\alpha H_LM(|x_j|)\longrightarrow0,
\]
a contradiction.

Define $V_j(x)=R_j^\alpha u(R_jx)$. Passing to a subsequence,
the annular estimates in Section~\ref{sec:prelim} give
\[
 V_j\longrightarrow V\quad\text{in }
 C^1_{\loc}(\R^n\setminus\{0\}),
\]
where $V$ solves $-\Delta_pV=V^q$ and satisfies
\[
 0\le V(x)\le C_0|x|^{-\alpha}\quad(x\ne0),
 \qquad\max_{|x|=1}V(x)=\varepsilon.
\]
By \eqref{eq:rem-crossing} and \eqref{eq:rem-crossing-separation},
\[
 V(x)\le\varepsilon|x|^{-\alpha}\qquad(0<|x|<1).
\]
Lemma~\ref{lem:rem-smallness} extends $V$ to a nonnegative $C^1$
weak solution of $-\Delta_pV=V^q$ in $\R^n$.
Under \eqref{eq:range}, the Liouville theorem of
Serrin--Zou~\cite[Corollary~II(ii)]{SerrinZou2002} gives $V\equiv0$,
contradicting $\max_{|x|=1}V(x)=\varepsilon>0$.
This proves the claim. Estimate \eqref{eq:rem-small-hypothesis}
therefore holds in a sufficiently small ball, and Lemma~\ref{lem:rem-smallness} removes the singularity of $u$.

Conversely, a removable solution is locally bounded, and hence
$|x|^\alpha u(x)\to0$. Consequently, a nonremovable solution has
$\liminf_{x\to0}|x|^\alpha u(x)>0$.
Combining this with \eqref{eq:upper} proves \eqref{eq:two-sided-local}.
\end{proof}

\section{A continuous Steiner deformation}\label{sec:steiner}
In this section, we construct a deformation adapted to a solution that
is unbounded at the origin. We follow the interval construction of
Brock~\cite{Brock1995,Brock2000}, but allow the velocity of an interval
to vanish when its center is close to the origin. This keeps the high
superlevel sets fixed. We establish both an energy inequality and a
strict decrease criterion for this deformation.

Fix a unit vector $e\in\Sph^{n-1}$ and write $x=y+ze$, where
$y\in e^\perp$ and $z\in\R$. Let $\varepsilon>0$, and choose a bounded
odd function $f\in C^\infty(\R)$ such that
\begin{equation}\label{eq:velocity-choice}
 f(c)=0\quad\text{if }|c|\le\varepsilon,
 \qquad f'(c)>0\quad\text{if }|c|>\varepsilon.
\end{equation}
Denote by $F_t$ the flow generated by $\dot c=-f(c)$,
i.e., $\partial_t F_t(c)=-f(F_t(c))$ with $F_0(c)=c$,
and put $M_f=\|f\|_\infty$.
For any interval $I=(c-\ell,c+\ell)$ present at time
$t_0\ge0$, we denote by $I(t)$ its position at time
$t\ge t_0$, before it meets another interval, and define
\[
I(t)=\bigl(F_{t-t_0}(c)-\ell,\,
           F_{t-t_0}(c)+\ell\bigr).
\]
We first define an auxiliary transformation $T_t$
on finite unions of bounded open intervals,
\[
A=\bigcup_{k=1}^N I_k,
\qquad I_k=(c_k-\ell_k,c_k+\ell_k).
\]
We take the $I_k$ to be the disjoint component intervals of $A$.
Before $I_k$ meets another interval, its position at time $t$ is
\[
I_k(t)=
\bigl(F_t(c_k)-\ell_k,\,F_t(c_k)+\ell_k\bigr).
\]
If, at some time $t_0\ge0$, two current intervals are
\[
(c-l,b)\qquad\text{and}\qquad(b,c+l),\qquad a<b<d,
\]
we replace them by the single interval $(c-l,c+l)$, and then move the new interval according to
\[
\bigl(F_{t-t_0}(c)-\ell,\,
      F_{t-t_0}(c)+\ell\bigr),
\qquad t\ge t_0,
\]
until it meets another interval.
If several consecutive intervals have matching endpoints
at the same time, we replace them by the single interval
between their two outermost endpoints.
We denote the resulting union of intervals by $T_tA$.

We next define $S_t$ on bounded open sets
$A\subset\mathbb R$.
Choose finite unions of bounded open intervals $A_j$ such that
\[
A_j\subset A_{j+1},\qquad
\overline{A_j}\Subset A,\qquad
\bigcup_{j=1}^{\infty}A_j=A,
\]
and set
\[
S_tA:=\bigcup_{j=1}^{\infty}T_tA_j.
\]

We now extend $S_t$ to nonnegative Lipschitz functions $v:\mathbb R^n\to[0,\infty)$ satisfying $v(x)\to0$ as $|x|\to\infty$.
By an abuse of notation, we also denote the induced transformation on functions by $S_t$.
For $y\in e^\perp$ and $s>0$, set
\[
A_{y,s}:=\{z\in\mathbb R:v(y+ze)>s\}.
\]
We define
\[
(S_tv)(y+ze)
:=
\sup\Bigl(
\{r\in\mathbb Q:\ r>0,\ z\in S_t A_{y,r}\}
\cup\{0\}
\Bigr).
\]
The following lemma establishes the properties
of this construction.

\begin{lemma}\label{lem:frozen-steiner}
For every bounded open set $A\subset\mathbb R$, the definition of $S_tA$
is independent of the chosen approximating sequence $(A_j)$.
The construction above defines an order-preserving,
equimeasurable semigroup $\{S_t\}_{t\ge0}$ on
nonnegative Lipschitz functions tending to zero
at infinity, with the following properties.
\begin{enumerate}
\item
For every $t\ge0$,
\begin{equation}\label{eq:frozen-lipschitz}
\operatorname{Lip}(S_tv)\le\operatorname{Lip}(v),
\end{equation}
and
\begin{equation}\label{eq:frozen-displacement}
\inf_{|h|\le M_ft}v(x+he)
\le S_tv(x)
\le
\sup_{|h|\le M_ft}v(x+he)
\end{equation}
for all $x\in\mathbb R^n$.

 \item If $\{v>T\}\subset\{|x\cdot e|<\varepsilon\}$, then
 \begin{equation}\label{eq:frozen-high-levels}
 \{S_tv>s\}=\{v>s\}\quad(s\ge T),
 \qquad S_tv=v\quad\text{on }\{v>T\}.
 \end{equation}
\end{enumerate}
\end{lemma}

\begin{proof}

We first show that $T_t$ preserves the total length
of finite unions of intervals and is order-preserving.
For any interval $I=(c-\ell,c+\ell)$, we have
\(
|I(t)|=2\ell=|I|
\)
before it meets another interval.
If two intervals $(a,b)$ and $(b,d)$ meet, they are
replaced by $(a,d)$ according to the construction.
Since
\(
|(a,d)|=|(a,b)|+|(b,d)|,
\)
the total length remains unchanged.
The resulting interval continues to move according to
the same rule until it meets another interval.
Thus $T_t$ preserves the total length of finite unions
of intervals.

Differentiating
$\partial_tF_t(c)=-f(F_t(c))$
with respect to $c$, we obtain
\[
\partial_tF_t'(c)
=-f'(F_t(c))F_t'(c),
\qquad F_0'(c)=1.
\]
Consequently, since $f'\ge0$,
\begin{equation}\label{eq:frozen-flow-derivative}
0<F_t'(c)
=\exp\left(-\int_0^t f'(F_\tau(c))\,\dd\tau\right)
\le1.
\end{equation}
For any $A\subset\mathbb R$ and $h\ge0$, write
\[
A+[-h,h]:=\{a+s:a\in A,\ |s|\le h\}.
\]
First consider two intervals
\[
A=(c_A-\ell_A,c_A+\ell_A),\qquad
B=(c_B-\ell_B,c_B+\ell_B).
\]
The condition $A+[-h,h]\subset B$ is equivalent to
\[
|c_A-c_B|\le\ell_B-\ell_A-h.
\]
By \eqref{eq:frozen-flow-derivative} and the mean value theorem,
\[
|F_t(c_A)-F_t(c_B)|
\le |c_A-c_B|
\le\ell_B-\ell_A-h.
\]
Thus $T_tA+[-h,h]\subset T_tB$ as long as
neither interval meets another interval.
For finite unions, this argument holds for each component of
$A$ and the component of $B$.
Induction over the finitely many mergers proves the inclusion
for all $t\ge0$.
Consequently, for any such $A,B$ and $h\ge0$,
\begin{equation}\label{eq:finite-buffer-order}
A+[-h,h]\subset B \quad\Longrightarrow\quad T_tA+[-h,h]\subset T_tB.
\end{equation}
Taking $h=0$ gives the order-preserving property of $T_t$.

We now verify the definition of $S_t$ on bounded open sets.
Let $A_j$ and $\widetilde A_k$ be two sequences used
in the definition for the same set $A$.
For each $j$, compactness gives
$\overline{A_j}\subset\widetilde A_k$ for some $k$.
By the order-preserving property of $T_t$, we have
\(
T_tA_j\subset T_t\widetilde A_k.
\)
Thus
\[
\bigcup_{j=1}^\infty T_tA_j
\subset
\bigcup_{k=1}^\infty T_t\widetilde A_k.
\]
Interchanging the two sequences shows that the definition of $S_tA$ is independent of the sequence.
Moreover, since $T_tA_j$ increases with $j$,
\[
|S_tA|
=\lim_{j\to\infty}|T_tA_j|
=\lim_{j\to\infty}|A_j|
=|A|.
\]
Let $A,B\subset\mathbb R$ be bounded open sets
satisfying $A+[-h,h]\subset B$, where $h\ge0$.
Choose approximating sequences $A_j$ and $B_k$
as in the definition.
For each $j$, the compact set
$\overline{A_j}+[-h,h]$ is contained in $B$,
and hence in some $B_k$.
The result for finite unions of intervals gives
\[
T_tA_j+[-h,h]\subset T_tB_k\subset S_tB.
\]
Taking the union over $j$, we obtain
\begin{equation}\label{eq:frozen-buffer-order}
A+[-h,h]\subset B
\quad\Longrightarrow\quad
S_tA+[-h,h]\subset S_tB.
\end{equation}
In particular, $S_t$ is order-preserving on bounded
open sets.

For $v$ as in the lemma, fix $y\in e^\perp$ and $s>0$.
By the definition of $S_t$ on bounded open sets
and its order-preserving property,
\[
S_tA_{y,s}
\subset
\bigcup_{\substack{r\in\mathbb Q\\r>s}}S_tA_{y,r}.
\]
The reverse inclusion follows from
$A_{y,r}\subset A_{y,s}$ whenever $r>s$.
Consequently,
\[
S_tA_{y,s}
=
\bigcup_{\substack{r\in\mathbb Q\\r>s}}S_tA_{y,r}.
\]
Thus
\begin{equation}\label{eq:frozen-superlevel-sets}
\{z\in\mathbb R:(S_tv)(y+ze)>s\}=S_tA_{y,s}
=S_t\{z:v(y+ze)>s\}.
\end{equation}
To prove that $\{S_t\}_{t\ge0}$ is a semigroup, fix a bounded open set $A$ and an approximating sequence $(A_j)$ as defined above.
For each $j$, choose $k$ and $h>0$ such that $\overline{A_j}+[-h,h]\subset A_k$.
By \eqref{eq:finite-buffer-order}, we obtain
\(
T_\tau A_j+[-h,h]\subset T_\tau A_k.
\)
Thus $\overline{T_\tau A_j}\Subset T_\tau A_k\subset S_\tau A$, and $(T_\tau A_j)_j$ is again an approximating sequence for $S_\tau A$.
Since $F_{t+\tau}=F_t\circ F_\tau$, the identity
$T_{t+\tau}A=T_t(T_\tau A)$ holds before
any two intervals meet.
When two intervals meet, they are replaced by a
single interval, which follows the same rule
starting from its position at that time.
Thus the identity remains valid after each such
event, and hence
\[
T_{t+\tau}A=T_t(T_\tau A)
\]
for every finite union of bounded open intervals $A$.
Since $(T_\tau A_j)_j$ is an approximating sequence
for $S_\tau A$, we obtain
\[
S_t(S_\tau A)=\bigcup_{j=1}^\infty T_t(T_\tau A_j)=\bigcup_{j=1}^\infty T_{t+\tau}A_j=S_{t+\tau}A.
\]
Also, $S_0A=A$ by the definition of $S_t$.
Using \eqref{eq:frozen-superlevel-sets}, we obtain
\[
S_0v=v,\qquad S_{t+\tau}v=S_t(S_\tau v)\qquad(t,\tau\ge0).
\]
Since $|f|\le M_f$, we have
\[
|F_t(c)-c|
\le \int_0^t |f(F_\tau(c))|\,\dd\tau
\le M_ft.
\]
For an interval $I=(c-\ell,c+\ell)$, before it meets
another interval,
\[
I(t)=(F_t(c)-\ell,F_t(c)+\ell).
\]
Thus both endpoints move with speed at most $M_f$. Thus \[
I(t)\subset I+[-M_ft,M_ft].
\]
When two intervals meet, they are replaced by a single
interval whose outer endpoints coincide with those of
the original intervals at that time. Hence the endpoints
do not jump and continue to move with speed at most $M_f$.
Thus
\begin{equation}\label{eq:frozen-set-bound}
   S_tA\subset A+[-M_f t,M_f t].
\end{equation}
Define
\[
A\ominus[-h,h]:=\{z\in\mathbb R:z+[-h,h]\subset A\}.
\]
Take $z\in A\ominus[-M_f t,M_f t]$.
Since $[z-M_f t,z+M_f t]$ is compactly contained in the
open set $A$, there exists $\sigma>0$ such that
$J=(z-\sigma,z+\sigma)$ satisfies
\(
J+[-M_f t,M_f t]\subset A.
\)
By \eqref{eq:frozen-buffer-order},
\[
S_tJ+[-M_f t,M_f t]\subset S_tA.
\]
As $J$ is a single interval,
\[
S_tJ=(F_t(z)-\sigma,F_t(z)+\sigma).
\]
Together with $|F_t(z)-z|\le M_f t$, this implies
\[
z\in S_tJ+[-M_f t,M_f t]\subset S_tA.
\]
Hence
\begin{equation}\label{eq:frozen-inner-bound}
A\ominus[-M_f t,M_f t]\subset S_tA.
\end{equation}
We now apply \eqref{eq:frozen-set-bound} and \eqref{eq:frozen-inner-bound} to $A_{y,s}$.
Fix $x=y+ze$.
If $(S_tv)(x)>s$, then $z\in S_tA_{y,s}$.
By \eqref{eq:frozen-set-bound}, we have
\[
S_tA_{y,s}\subset A_{y,s}+[-M_f t,M_f t].
\]
Then there exists $|h|\le M_f t$ such that
$v(x+he)>s$. Hence
\[
(S_tv)(x)\le\sup_{|h|\le M_f t}v(x+he).
\]
Conversely, if
\[
\inf_{|h|\le M_f t}v(x+he)>s,
\]
then $[z-M_f t,z+M_f t]\subset A_{y,s}$, so
$z\in A_{y,s}\ominus[-M_f t,M_f t]\subset S_tA_{y,s}$.
Consequently, $(S_tv)(x)>s$.
Letting $s$ approach the infimum from below, we obtain
\[
\inf_{|h|\le M_f t}v(x+he)\le(S_tv)(x).
\]
Combining the two inequalities proves \eqref{eq:frozen-displacement}.

To prove \eqref{eq:frozen-lipschitz}, let $L$ be a Lipschitz constant of
$v$. For $y_1,y_2\in e^\perp$ and $h\ge0$,
\[
 \{v(y_1+\cdot e)>s+L\sqrt{|y_1-y_2|^2+h^2}\}+[-h,h]
 \subset\{v(y_2+\cdot e)>s\}.
\]
Indeed, let $z$ belong to the set on the left.
By the definition of the sum of sets, there exists
$z_0\in\mathbb R$ such that $|z-z_0|\le h$ and
\[
v(y_1+z_0e)>
s+L\sqrt{|y_1-y_2|^2+h^2}.
\]
Since $y_1-y_2\perp e$, the Lipschitz continuity
of $v$ gives
\[
\begin{aligned}
v(y_2+ze)
&\ge v(y_1+z_0e)
-L|(y_2+ze)-(y_1+z_0e)|\\
&\ge v(y_1+z_0e)
-L\sqrt{|y_1-y_2|^2+h^2}\\
&>s.
\end{aligned}
\]
This proves the claimed inclusion.
Applying \eqref{eq:frozen-buffer-order}, we obtain
\[
S_tA_{y_1,s+LD}+[-h,h]
\subset S_tA_{y_2,s},
\qquad
D=\sqrt{|y_1-y_2|^2+h^2}.
\]
For any $z_1,z_2\in\mathbb R$, take
$h=|z_1-z_2|$. By the definition of $S_tv$,
the preceding inclusion implies
\begin{equation}\label{eq:frozen-level-implication}
(S_tv)(y_1+z_1e)>s+LD\quad\Longrightarrow\quad(S_tv)(y_2+z_2e)>s.
\end{equation}
Suppose that
\[
(S_tv)(y_1+z_1e)>(S_tv)(y_2+z_2e)+LD.
\]
Since $S_tv\ge0$, we may choose $s>0$ such that
\[
(S_tv)(y_2+z_2e)<s<(S_tv)(y_1+z_1e)-LD.
\]
By \eqref{eq:frozen-level-implication}, we obtain
$(S_tv)(y_2+z_2e)>s$, a contradiction. Hence
\[
(S_tv)(y_1+z_1e)
\le (S_tv)(y_2+z_2e)+LD.
\]
Interchanging the two points gives
\[
\begin{aligned}
&|(S_tv)(y_1+z_1e)-(S_tv)(y_2+z_2e)|\le
L\sqrt{|y_1-y_2|^2+|z_1-z_2|^2}.
\end{aligned}
\]
Hence
\[
\operatorname{Lip}(S_tv)\le\operatorname{Lip}(v).
\]
In particular, $S_tv$ is Lipschitz continuous and hence measurable.

If $\{v>T\}\subset\{|x\cdot e|<\varepsilon\}$, then for every $s\ge T$, each interval $I=(c-\ell,c+\ell)$ contained in $A_{y,s}$ satisfies $|c|<\varepsilon$.
Since $f(c)=0$, the constant function $g(t)\equiv c$ solves the initial value problem for $F_t(c)$.
By uniqueness, $F_t(c)=c$ for all $t\ge0$. Thus
\[
S_tA_{y,s}=A_{y,s},\qquad s\ge T.
\]
Consequently,
\[
\{S_tv>s\}=\{v>s\}
\qquad(s\ge T),
\]
which also implies $S_tv=v$ on $\{v>T\}$.
This proves \eqref{eq:frozen-high-levels}.

\end{proof}

\begin{lemma}\label{lem:steiner-energy}
Let $v\ge0$ be a Lipschitz function tending to zero at infinity,
and suppose that $\int_{\R^n}|\nabla v|^p\dd x<\infty$.
Then, for every $t\ge0$,
 \begin{equation}\label{eq:steiner-energy}
 \int_{\R^n}|\nabla S_tv|^p\dd x
 \le\int_{\R^n}|\nabla v|^p\dd x.
 \end{equation}
\end{lemma}

\begin{proof}
We prove a stronger version of \eqref{eq:steiner-energy}.
For a Borel set $D\subset e^\perp\times(0,\infty)$, define
\[
 \mu_v(D)=\int_{\{(y,z):(y,v(y+ze))\in D\}}
             |\nabla v(y+ze)|^p\,\dd y\,\dd z.
\]
We will show that
\begin{equation}\label{eq:frozen-weighted-polya}
 \mu_{S_tv}(D)\le\mu_v(D).
\end{equation}
Let $w\ge0$ be compactly supported and Lipschitz, and put $p'=p/(p-1)$. For $d>0$, define
\[
 H_dw(x)=\sup_{\xi\in\R^n}
 \left\{w(x+\xi)-\frac{|\xi|^{p'}}{p'd^{p'-1}}\right\}.
\]
We claim that
\begin{equation}\label{eq:frozen-envelope}
 H_d(S_tw)\le S_t(H_dw).
\end{equation}
Fix $y,y'\in e^\perp$ and $r\ge0$.
Let $z$ belong to the set on the left.
Then there exists $z_0\in\mathbb R$ such that
\[
|z-z_0|\le r,
\qquad
w(y'+z_0e)>a+k.
\]
Taking
$\xi=(y'-y)+(z_0-z)e$,
we have
\[
|\xi|^2
=|y-y'|^2+|z_0-z|^2
\le |y-y'|^2+r^2.
\]
By the definition of $H_d$,
\[
\begin{aligned}
H_dw(y+ze)
\ge w(y'+z_0e)
-\frac{|\xi|^{p'}}{p'd^{p'-1}}\ge w(y'+z_0e)-k>a.
\end{aligned}
\]
This proves that
\[
 \left\{w(y'+\cdot e)>a+\frac{(|y-y'|^2+r^2)^{p'/2}}{p'd^{p'-1}}\right\}+[-r,r]
 \subset\{H_dw(y+\cdot e)>a\}\qquad(a>0).
\]
Applying \eqref{eq:frozen-buffer-order}, we obtain
\[
\begin{aligned}
\{z:(S_tw)(y'+ze)>a+\frac{(|y-y'|^2+r^2)^{p'/2}}{p'd^{p'-1}}\}+[-r,r]
\subset
\{z:S_t(H_dw)(y+ze)>a\}.
\end{aligned}
\]
Now fix $x=y+ze$ and suppose that
$H_d(S_tw)(x)>a$ for some $a>0$.
By the definition of $H_d$, there exists
$\xi\in\mathbb R^n$ such that
\[
(S_tw)(x+\xi)>
a+\frac{|\xi|^{p'}}{p'd^{p'-1}}.
\]
Write $x+\xi=y'+z_0e$ and choose
$r=|z-z_0|$. Then
\[
(S_tw)(y'+z_0e)>a+\frac{(|y-y'|^2+r^2)^{p'/2}}{p'd^{p'-1}},
\qquad |z-z_0|=r.
\]
The preceding inclusion gives
$S_t(H_dw)(x)>a$.
Since $a>0$ is arbitrary, we obtain \eqref{eq:frozen-envelope}.

Let $L=\operatorname{Lip}(w)$.
 By the Lipschitz continuity of $w$, we have
\[
w(x+\xi)\le w(x)+L|\xi|.
\]
If $|\xi|>(p'L)^{p-1}d$, then
\[
\frac{|\xi|^{p'}}{p'd^{p'-1}}>L|\xi|,
\]
and hence
\[
w(x+\xi)-\frac{|\xi|^{p'}}{p'd^{p'-1}}
<w(x).
\]
On the other hand, taking $\xi=0$ in the definition
of $H_dw$ gives $H_dw(x)\ge w(x)$.
Therefore, the supremum can be restricted to
$|\xi|\le(p'L)^{p-1}d$, namely,
\begin{equation}\label{eq:frozen-envelope-restriction}
H_dw(x)=\sup_{|\xi|\le(p'L)^{p-1}d}\left\{w(x+\xi)-\frac{|\xi|^{p'}}{p'd^{p'-1}}\right\}.
\end{equation}
Setting $\xi=d\zeta$, we obtain
\[
\frac{H_dw(x)-w(x)}{d}
=\sup_{\zeta\in\R^n}
\left\{
\frac{w(x+d\zeta)-w(x)}{d}
-\frac{|\zeta|^{p'}}{p'}
\right\}.
\]
Consequently,
\begin{equation}\label{eq:frozen-envelope-bound}
0\le\frac{H_dw(x)-w(x)}{d}
\le\sup_{r\ge0}
\left\{Lr-\frac{r^{p'}}{p'}\right\}
=\frac{L^p}{p}.
\end{equation}
Now let $x$ be a differentiability point of $w$. By \eqref{eq:frozen-envelope-restriction},
\[
\frac{H_dw(x)-w(x)}{d}
=\sup_{|\zeta|\le (p'L)^{p-1}}
\left\{
\frac{w(x+d\zeta)-w(x)}{d}
-\frac{|\zeta|^{p'}}{p'}
\right\}.
\]
Since $w$ is differentiable at $x$, for every
$\varepsilon>0$, there exists $\delta>0$ such that
\[
|w(x+h)-w(x)-\nabla w(x)\cdot h|
\le \varepsilon |h|
\qquad (|h|<\delta).
\]
If $L>0$ and $0<d<{\delta}/{(p'L)^{p-1}}$, then for every
$|\zeta|\le (p'L)^{p-1}$, we have $|d\zeta|<\delta$ and hence
\[
\begin{aligned}
\left|
\frac{w(x+d\zeta)-w(x)}{d}
-\nabla w(x)\cdot\zeta
\right|&=
\frac{1}{d}
\left|
w(x+d\zeta)-w(x)
-\nabla w(x)\cdot(d\zeta)
\right|\\
&\le \varepsilon|\zeta|
\le \varepsilon (p'L)^{p-1}.
\end{aligned}
\]
Taking the supremum over $|\zeta|\le (p'L)^{p-1}$ gives
\[
\sup_{|\zeta|\le (p'L)^{p-1}}
\left|
\frac{w(x+d\zeta)-w(x)}{d}
-\nabla w(x)\cdot\zeta
\right|
\le \varepsilon (p'L)^{p-1}.
\]
Since $\varepsilon>0$ is arbitrary, we conclude that
\[
\sup_{|\zeta|\le (p'L)^{p-1}}
\left|
\frac{w(x+d\zeta)-w(x)}{d}
-\nabla w(x)\cdot\zeta
\right|
\longrightarrow 0
\qquad(d\downarrow0).
\]
Hence we may pass to the limit under the supremum:
\[
\begin{aligned}
\lim_{d\downarrow0}\frac{H_dw(x)-w(x)}{d}
&=\sup_{|\zeta|\le (p'L)^{p-1}}
\left\{
\nabla w(x)\cdot\zeta
-\frac{|\zeta|^{p'}}{p'}
\right\}.
\end{aligned}
\]
The last supremum is attained at $\zeta=|\nabla w(x)|^{p-2}\nabla w(x)$ if $\nabla w(x)\ne0$, and at $\zeta=0$ otherwise.
This point lies in the ball $|\zeta|\le (p'L)^{p-1}$, since $|\nabla w(x)|\le L$.
Thus we conclude that
\begin{equation}\label{eq:frozen-envelope-limit}
\lim_{d\downarrow0}\frac{H_dw(x)-w(x)}{d}=\frac{|\nabla w(x)|^p}{p}\quad\text{a.e. in }\mathbb R^n.
\end{equation}
By \eqref{eq:frozen-lipschitz}, equations~\eqref{eq:frozen-envelope-bound} and~\eqref{eq:frozen-envelope-limit} also hold for $S_tw$.
By fiberwise equimeasurability and the layer-cake formula, for every nonnegative integrable function $v$ considered here,
\[
\begin{aligned}
\int_{\mathbb R}(S_tv)(y+ze)\,dz
&=\int_0^\infty
|\{z:(S_tv)(y+ze)>s\}|\,ds=\int_0^\infty
|\{z:v(y+ze)>s\}|\,ds\\
&=\int_{\mathbb R}v(y+ze)\,dz.
\end{aligned}
\]
Thus for any measurable $Y\subset e^\perp$,
\[
\int_{Y\times\mathbb R}S_t(H_dw)
=
\int_{Y\times\mathbb R}H_dw,\qquad
\int_{Y\times\mathbb R}S_tw
=
\int_{Y\times\mathbb R}w.
\]
Integrate \eqref{eq:frozen-envelope} over $Y\times\mathbb R$ after subtracting $S_tw$, we obtain
\[
\begin{aligned}
\int_{Y\times\mathbb R}
\bigl(H_d(S_tw)-S_tw\bigr)\,\dd y\,\dd z  &\le
\int_{Y\times\mathbb R}
\bigl(S_t(H_dw)-S_tw\bigr)\,\dd y\,\dd z\\
&=
\int_{Y\times\mathbb R}
(H_dw-w)\,\dd y\,\dd z.
\end{aligned}
\]
Choose $R>0$ such that $w=0$ outside $B_R$.
By \eqref{eq:frozen-envelope-restriction}, for $0<d\le1$ we have
\[
|\xi|\le (p'L)^{p-1}d\le (p'L)^{p-1}.
\]
Thus, if $|x|>R+(p'L)^{p-1}$, then
\[
|x+\xi|\ge |x|-|\xi|>R.
\]
Since $w=0$ outside $B_R$, it follows that
$w(x+\xi)=0$ for all such $\xi$. Consequently,
\[
\begin{aligned}
H_dw(x)
=\sup_{|\xi|\le (p'L)^{p-1}d}
\left\{
-\frac{|\xi|^{p'}}{p'd^{p'-1}}
\right\}=0,
\end{aligned}
\]
where the last equality follows by taking $\xi=0$.
Therefore,
\[
H_dw=0
\qquad\text{outside }B_{R+(p'L)^{p-1}}.
\]
By \eqref{eq:frozen-displacement}, $S_tw=0$ outside $B_{R+M_ft}$. The same argument gives
\[
H_d(S_tw)=0
\quad\text{outside }B_{R+M_ft+(p'L)^{p-1}}
\]
for every $0<d\le1$.
Thus both difference quotients
\[
\frac{H_dw-w}{d},
\qquad
\frac{H_d(S_tw)-S_tw}{d}
\]
are supported in the same compact set
$\overline{B_{R+M_ft+(p'L)^{p-1}}}$, independently of $d$.
By \eqref{eq:frozen-envelope-bound}, they are bounded
by the integrable function
\[
\frac{L^p}{p}\,
\mathbf{1}_{B_{R+M_ft+(p'L)^{p-1}}}.
\]
Applying \eqref{eq:frozen-envelope-limit} and dominated convergence, we obtain
\begin{equation}\label{eq:frozen-horizontal}
 \int_{Y\times\R}|\nabla S_tw|^p
 \le\int_{Y\times\R}|\nabla w|^p.
\end{equation}
Now let $v$ satisfy the hypotheses of this lemma.
Given $0<l<h$, let $w=\min\{(v-l)_+,h-l\}$.
Since $v(x)\to0$ as $|x|\to\infty$, the function $w$ is compactly supported.
Moreover, the definition of $S_t$ gives
\[
S_tw=\min\{(S_tv-l)_+,h-l\},
\]
since for every $0<s<h-l$,
\[
\{w>s\}=\{v>l+s\}.
\]
Since
\[
\nabla w=\mathbf{1}_{\{l<v<h\}}\nabla v\qquad
\nabla S_tw
=\mathbf{1}_{\{l<S_tv<h\}}\nabla S_tv
\quad\text{a.e.}
\]
we have
\begin{equation}\label{eq:frozen-truncation-measures}
\begin{aligned}
\int_{Y\times\mathbb R}|\nabla w|^p\,dx
&=\int_{Y\times\mathbb R}
\mathbf1_{\{l<v<h\}}|\nabla v|^p\,dx=\mu_v(Y\times(l,h)),\\
\int_{Y\times\mathbb R}|\nabla S_tw|^p\,\dd x
&=\int_{Y\times\mathbb R}\mathbf1_{\{l<S_tv<h\}}|\nabla S_tv|^p\,\dd x
=\mu_{S_tv}(Y\times(l,h)).
\end{aligned}
\end{equation}
Hence \eqref{eq:frozen-horizontal} implies
\[
\mu_{S_tv}(Y\times(l,h))\le\mu_v(Y\times(l,h)).
\]
Since $v$ and $S_tv$ are Lipschitz continuous, their
gradients vanish almost everywhere on each level set.
Hence, for every $s>0$ and measurable $Y\subset e^\perp$,
\[
\mu_v(Y\times\{s\})
=\mu_{S_tv}(Y\times\{s\})=0.
\]
Consequently, \eqref{eq:frozen-truncation-measures} also holds
for $Y\times(l,h]$.
Let
\[
D_m=\{y\in e^\perp:|y|<m\}\times(1/m,m).
\]
Since $D_m$ increases to
$\Omega=e^\perp\times(0,\infty)$, we have
\[
\begin{aligned}
\mu_{S_tv}(\Omega)
&=\lim_{m\to\infty}\mu_{S_tv}(D_m)\le\lim_{m\to\infty}\mu_v(D_m)
=\mu_v(\Omega)\le\int_{\R^n}|\nabla v|^p\,\dd x<\infty.
\end{aligned}
\]
Thus \(\mu_{S_tv}\) and \(\mu_v\) are finite measures.
By additivity, the inequality
$\mu_{S_tv}(D)\le\mu_v(D)$ holds for finite
unions of disjoint rectangles of the form
$Y\times(l,h]$. By taking limits, it also holds
when $l=0$ or $h=\infty$.
Now let
\[
\mathcal M=
\{D\subset\Omega:D\text{ is Borel and }
\mu_{S_tv}(D)\le\mu_v(D)\}.
\]
Since both measures are finite, $\mathcal M$
is closed under increasing unions and decreasing
intersections, by continuity of measures.
Moreover, $\mathcal M$ contains the algebra
of finite unions of the above rectangles.
Since this algebra generates the Borel
$\sigma$-algebra of $\Omega$, the monotone class
theorem gives
\[
\mu_{S_tv}(D)\le\mu_v(D)
\]
for every Borel set $D\subset\Omega$.
This proves \eqref{eq:frozen-weighted-polya}.

Finally, since the gradients of $v$ and $S_tv$
vanish almost everywhere on their zero sets,
\[
\mu_v(\Omega)=\int_{\R^n}|\nabla v|^p\,\dd x,
\qquad
\mu_{S_tv}(\Omega)
=\int_{\R^n}|\nabla S_tv|^p\,\dd x.
\]
Taking $D=\Omega$ in
\eqref{eq:frozen-weighted-polya} yields
\eqref{eq:steiner-energy}.

\end{proof}

The next lemma gives the strict energy decrease needed to obtain
symmetry. The derivatives below are taken with respect to both $y$
and the level $s$.

\begin{lemma}\label{lem:strict-steiner}
Let $v$ satisfy the hypotheses of Lemma~\ref{lem:steiner-energy}.
Suppose there exist $C^1$ functions
$a_-,a_+:\mathcal V\to\mathbb R$, where
$\mathcal V\subset e^\perp\times(0,\infty)$ is open,
such that
\[
a_-(y,s)<a_+(y,s),
\]
and $(a_-(y,s),a_+(y,s))$ is a connected component
of $\{z\in\mathbb R:v(y+ze)>s\}$
for every $(y,s)\in\mathcal V$.
Assume that $v$ is $C^1$ near these endpoints and
\[
 \partial_ev(y+a_-e)>0,\qquad
 \partial_ev(y+a_+e)<0.
\]
Set $c=(a_-+a_+)/2$. If at some point of $\mathcal V$,
\begin{equation}\label{eq:frozen-nontrivial-center}
 f'(c)>0,\qquad D_{y,s}c\ne0,
\end{equation}
then there are $b,t_0>0$ such that
\begin{equation}\label{eq:frozen-strict-deficit}
 \int_{\R^n}|\nabla S_tv|^p\dd x
 \le\int_{\R^n}|\nabla v|^p\dd x-bt
 \qquad(0<t<t_0).
\end{equation}
\end{lemma}

\begin{proof}
Let $(y_0,s_0)\in\mathcal V$ satisfy
\[
f'(c(y_0,s_0))>0,
\qquad
D_{y,s}c(y_0,s_0)\ne0.
\]
By continuity, choose a sufficiently small open cube $Q\subset e^\perp$ containing $y_0$ and $0<l<s_0<h$ such that
\[
D=\overline Q\times[l,h]\Subset\mathcal V,
\]
and
\[
f'(c)>0,\qquad D_{y,s}c\ne0
\quad\text{on }D.
\]
By continuity and the compactness of $D$, there exists
$\kappa>0$ such that
\[
\partial_e v(y+a_-(y,s)e)\ge\kappa,
\qquad
\partial_e v(y+a_+(y,s)e)\le-\kappa
\]
for all $(y,s)\in D$.
We next choose two fixed open intervals on either side of the interval $(a_-(y,s),a_+(y,s))$.
 Since
\[
\partial_e v(y_0+a_-(y_0,s_0)e)>0,
\qquad
\partial_e v(y_0+a_+(y_0,s_0)e)<0,
\]
we can choose
\[
z_-<a_-(y_0,s_0),\qquad
z_+>a_+(y_0,s_0)
\]
sufficiently close to the corresponding endpoints
such that
\[
v(y_0+z_-e)<s_0,\qquad
v(y_0+z_+e)<s_0.
\]
Choosing $l<s_0<h$ sufficiently close to $s_0$
and shrinking $Q$ if necessary, the continuity
of $v$ and $a_\pm$ allows us to find two fixed
open intervals $J_-$ and $J_+$ containing
$z_-$ and $z_+$, respectively, such that
\[
\sup J_-<a_-(y,s)<a_+(y,s)<\inf J_+
\]
for every $(y,s)\in D$, and
\[
v(y+ze)<l
\qquad
\text{for all }y\in\overline Q,\quad
z\in J_-\cup J_+.
\]
Between the endpoint neighborhoods, the compact interior of the chord stays above $h$.

Set $w=\min\{(v-l)_+,h-l\}$. Choose a Lipschitz cutoff $0\le\psi(z)\le1$
equal to one between the two exterior intervals and zero beyond them,
with its transitions inside those intervals. On $Q$, the supports
of $w_1=\psi w$ and $w_2=(1-\psi)w$ are separated by a fixed positive
distance. Since the interval speeds are bounded by $M_f$, their
transformed supports remain disjoint for all sufficiently small $t$.
The interval construction and compact exhaustion therefore give
\[
 S_tw=S_tw_1+S_tw_2\quad\text{on }Q\times\R.
\]
The Dirichlet energies are additive there before and after deformation,
and \eqref{eq:frozen-horizontal} controls the contribution from $w_2$.

It remains to compute the contribution from $w_1$. Write
$a_\pm=c\pm\ell$ and $A=-\ell_s>|c_s|$. Changing variables from $z$
to $s$ on the two inverse branches gives the energy density
\begin{equation}\label{eq:frozen-G}
 G(\xi,\eta)=
 \frac{(1+|\xi-\nabla_y\ell|^2)^{p/2}}{(A+\eta)^{p-1}}
 +\frac{(1+|\xi+\nabla_y\ell|^2)^{p/2}}{(A-\eta)^{p-1}},
 \qquad(\xi,\eta)=D_{y,s}c.
\end{equation}
The deformation replaces $c$ by $F_t(c)$ and leaves $\ell$ unchanged.
Thus $Dc$ is replaced by $F_t'(c)Dc$. The endpoint derivatives retain
their signs because $|F_t'(c)c_s|<A$.

For fixed $\nabla_y\ell$ and $A$, the function $G$ is even and
strictly convex on $|\eta|<A$. Indeed, $|X|^p/r^{p-1}$ is convex for $r>0$,
and equality in its convexity requires the ratios $X/r$ to agree.
Here the first component of $X$ is one, so equality forces both
arguments to agree. Since $Dc\ne0$,
\[
 \left.\frac{\mathrm d}{\mathrm dt}G(F_t'(c)Dc)\right|_{t=0}
 =-f'(c)DG(Dc)\cdot Dc<0.
\]
Compactness makes this decrease uniform on $D$ for small $t$.
Integration gives an energy loss of at least $bt$ for $w_1$ on
$Q\times\R$. Adding the estimate for $w_2$ gives the same loss
on the parameter set $Q\times(l,h)$. The two boundary levels have
zero energy mass, and \eqref{eq:frozen-weighted-polya} controls
the complement. This proves \eqref{eq:frozen-strict-deficit}.
\end{proof}

\section{Classification of the limiting solutions}\label{sec:rigidity}
We now classify positive solutions in the punctured whole space with
two-sided power bounds. This will identify every limit obtained by
rescaling a nonremovable solution near the origin.

\begin{proposition}\label{prop:global-rigidity}
Assume that \eqref{eq:range} holds. Let $U>0$ be a weak solution of
$-\Delta_pU=U^q$ in $\R^n\setminus\{0\}$. If
\begin{equation}\label{eq:global-bounds}
 c|x|^{-\alpha}\le U(x)\le C|x|^{-\alpha}
 \qquad(x\in\R^n\setminus\{0\})
\end{equation}
for some $c,C>0$, then $U(x)=\K|x|^{-\alpha}$.
\end{proposition}

\begin{proof}
The local regularity and scaling estimates in
Section~\ref{sec:prelim} give
\begin{equation}\label{eq:global-gradient}
 |\nabla U(x)|\le C|x|^{-\alpha-1}\qquad(x\ne0).
\end{equation}
We set $U(0)=+\infty$ when considering its superlevel sets.
For $M>0$, the truncation $V_M=\min\{U,M\}$ is bounded and globally
Lipschitz. It is constant near the origin and tends to zero at infinity.
Since $p(\alpha+1)=\alpha(q+1)=n+\delta>n$,
\begin{equation}\label{eq:cap-integrability}
 \int_{\R^n}\bigl(|\nabla V_M|^p+V_M^{q+1}\bigr)\dd x<\infty.
\end{equation}
We divide the proof into three steps.

\smallskip\noindent
\textit{Step 1. The energy variation is of order $O(t^2)$.}
Fix $e\in\Sph^{n-1}$ and $\varepsilon>0$. Choose $f$ as in
\eqref{eq:velocity-choice}, with $\|f\|_\infty\le1$, and let $S_t$
be the deformation in Lemma~\ref{lem:frozen-steiner}.
Take $C\varepsilon^{-\alpha}<m<M$, where $C$ is the upper constant
in \eqref{eq:global-bounds}. Then $\{V_M>m\}\subset B_\varepsilon$,
so \eqref{eq:frozen-high-levels} gives
\begin{equation}\label{eq:fixed-core}
 \phi_t:=S_tV_M-V_M=0\quad\text{on }\{U>m\}.
\end{equation}
In particular, $\phi_t$ vanishes in a neighborhood of the origin.
By \eqref{eq:frozen-displacement} and \eqref{eq:global-gradient},
\begin{equation}\label{eq:variation-tail}
 |\phi_t(x)|\le Ct|x|^{-\alpha-1}
 \qquad(|x|\ge2,\ 0<t<1/2).
\end{equation}
On a fixed compact set, $V_M$ has a positive minimum and
$|\phi_t|\le\operatorname{Lip}(V_M)t$. Together with the lower bound
in \eqref{eq:global-bounds}, this gives
$|\phi_t|\le C_{\varepsilon,M}tV_M$ on $\R^n$.

Put $F(s)=s^{q+1}/(q+1)$. Equimeasurability and
\eqref{eq:cap-integrability} imply
\begin{equation}\label{eq:potential-cancel}
 \int_{\R^n}\bigl[F(V_M+\phi_t)-F(V_M)\bigr]\dd x=0.
\end{equation}
For small $t$, apply $(1+h)^{q+1}=1+(q+1)h+O(h^2)$ with
$h=\phi_t/V_M$. Since $V_M^q\phi_t=U^q\phi_t$ almost everywhere,
we obtain
\begin{equation}\label{eq:taylor-cancel}
 \left|\int_{\R^n}U^q\phi_t\dd x\right|
 \le C_{\varepsilon,M}t^2.
\end{equation}
The relative bound on $\phi_t$ makes this estimate valid also for
$0<q<1$.

By Lemma~\ref{lem:steiner-energy}, $\phi_t$ has finite gradient
energy. Together with \eqref{eq:fixed-core}, this permits its use in
the weak equation after a cutoff at infinity. The cutoff error on
$B_{2R}\setminus B_R$ is bounded by
\[
 CtR^{n-p(\alpha+1)-1}=CtR^{-\delta-1}\longrightarrow0.
\]
On $\{U>m\}$, $\nabla\phi_t=0$ almost everywhere, and elsewhere
$\nabla V_M=\nabla U$. Writing
$E(v)=p^{-1}\int_{\R^n}|\nabla v|^p\dd x$, convexity gives
\begin{equation}\label{eq:energy-o-tau}
 0\ge E(S_tV_M)-E(V_M) \ge\int_{\R^n}|\nabla U|^{p-2}\nabla U\cdot\nabla\phi_t\dd x =\int_{\R^n}U^q\phi_t\dd x \ge-C_{\varepsilon,M}t^2.
\end{equation}

\smallskip\noindent
\textit{Step 2. $U$ is radially symmetric about the origin.}
Consider a component $(a_-(y,s),a_+(y,s))$ of a fiber superlevel
set of $U$, with regular endpoints. The implicit function theorem
gives $C^1$ endpoint graphs near $(y,s)$. These graphs enclose the
same component for nearby parameters: away from the endpoint
neighborhoods, the minimum of a continuous truncation of $U$ on the
compact interior of the chord is strictly greater than $s$.
This argument also applies when the chord passes through the origin.

Set $c(y,s)=(a_-+a_+)/2$. If $c\ne0$ and $Dc\ne0$ at some point,
choose $0<\varepsilon<|c|$ and take $m$ above both
$C\varepsilon^{-\alpha}$ and all levels in a small parameter
neighborhood. For $M>m$, the relevant chords of $U$ and $V_M$ agree.
Lemma~\ref{lem:strict-steiner} gives a linear energy decrease, which
contradicts \eqref{eq:energy-o-tau}. Hence $Dc=0$ wherever $c\ne0$.
By continuity, $Dc=0$ also where $c=0$; otherwise there would be
nearby points with both $c\ne0$ and $Dc\ne0$.
Differentiating $U(y+a_\pm e)=s$, we obtain
\begin{equation}\label{eq:local-reflection}
 \partial_eU(y+a_-e)=-\partial_eU(y+a_+e),
 \qquad
 \nabla_yU(y+a_-e)=\nabla_yU(y+a_+e).
\end{equation}

We also need these identities when only the entering endpoint is
known to be regular. Let $(a,b)$ be such a chord at height $s$.
By the one-dimensional Sard theorem, choose regular fiber values
$s_j\downarrow s$. Consider the component at height $s_j$ containing
a fixed interior point of $(a,b)$. Since $U>s$ on the chord, this
component eventually contains every compact subinterval of $(a,b)$.
Its endpoints therefore converge to $a,b$. Passing to the limit in
\eqref{eq:local-reflection} proves that the other endpoint is regular
and that both identities hold for the original chord. All endpoints
at finite levels are away from the singularity.

To apply the geometric local symmetry theorem of Brock, define
\[
 W=1-e^{-U}\quad\text{on }\R^n\setminus\{0\},
 \qquad W(0)=1.
\]
The bounds \eqref{eq:global-bounds}--\eqref{eq:global-gradient} give
$1-W=o(|x|)$ and
\[
 |\nabla W(x)|\le Ce^{-c|x|^{-\alpha}}|x|^{-\alpha-1}
 \longrightarrow0\quad\text{as }x\to0.
\]
Thus $W\in C^1(\R^n)$, $W>0$, and $W(x)\to0$ at infinity.
At paired endpoints of height $s$, its gradients are $e^{-s}$ times
those of $U$. Since $e$ was arbitrary, \eqref{eq:local-reflection}
gives local symmetry in every direction in the sense
of~\cite[Definitions~6.1--6.2 and Remark~6.1(a)]{Brock2000}.
By~\cite[Theorem~6.1]{Brock2000}, the set $\{\nabla W\ne0\}$ is a
disjoint union of open annuli on each of which $W$, and hence $U$,
is radial and strictly decreasing.

Fix one such annulus,
\begin{equation}\label{eq:radial-annulus}
 A=\{x:r<|x-z|<R\},\qquad
 U(x)=V(|x-z|),\qquad V'<0.
\end{equation}
The radial equation yields
\[
 H(\rho):=\rho^{n-1}(-V'(\rho))^{p-1}>0,
 \qquad H'(\rho)=\rho^{n-1}V(\rho)^q>0.
\]
Suppose $R<\infty$. At every point of the outer sphere distinct
from the origin, the gradient vanishes. Otherwise that point would
belong to another open annulus meeting $A$, contrary to disjointness.
Continuity would give $H(R^-)=0$, contradicting $H>0$ and $H'>0$.
Therefore each annulus is the exterior of a ball. Two such exteriors
intersect, so there is only one. There is at least one because $U$
is not constant.

If its inner radius $r$ were positive, then $\nabla U=0$ in
$B_r(z)\setminus\{0\}$. Testing the equation there with a nonzero
nonnegative smooth function would contradict $U^q>0$.
Hence $r=0$. Finally,
$A=\R^n\setminus\{z\}\subset\R^n\setminus\{0\}$ forces $z=0$.

\smallskip\noindent
\textit{Step 3. The radial solution is $\K r^{-\alpha}$.}
For the radial solution from Step~2, set
\[
 t=-\log r,\qquad v(t)=r^\alpha U(r),\qquad
 z(t)=v_t+\alpha v=-r^{\alpha+1}U'(r).
\]
By \eqref{eq:global-gradient} and $\beta>0$, the radial flux tends
to zero at the origin. Integrating the radial equation gives
\begin{equation}\label{eq:radial-flux}
 r^{n-1}(-U'(r))^{p-1}
 =\int_0^r\rho^{n-1}U(\rho)^q\dd\rho
 \asymp r^\beta.
\end{equation}
Here $\asymp$ means that the ratio is bounded above and below by
positive constants independent of $r$. Thus $v$ and $z$ remain
in compact subintervals of $(0,\infty)$ for all $t\in\R$, and
\begin{equation}\label{eq:radial-system}
 v_t=z-\alpha v,\qquad
 (z^{p-1})_t=\beta z^{p-1}-v^q.
\end{equation}
Define
\begin{equation}\label{eq:radial-energy}
 \mathcal E(v,z)=\frac{p-1}{p}z^p-\alpha vz^{p-1}
 +\frac{v^{q+1}}{q+1}+\frac{\delta\alpha^{p-1}}p v^p.
\end{equation}
Using \eqref{eq:radial-system} and $\delta=\alpha-\beta$, a direct
calculation gives
\begin{equation}\label{eq:radial-dissipation}
 \frac{\mathrm d\mathcal E}{\mathrm dt}
 =-\delta\bigl[z^{p-1}-(\alpha v)^{p-1}\bigr](z-\alpha v)\le0.
\end{equation}
Equality holds precisely when $z=\alpha v$.

The energy is bounded and nonincreasing, so it has limits as
$t\to\pm\infty$. Along either end, time translates of $(v,z)$
have subsequences converging on compact intervals to solutions of
\eqref{eq:radial-system}. The energy is constant on every such
limiting orbit. Therefore \eqref{eq:radial-dissipation} gives
$z=\alpha v$, and the system yields $v=\K$.
Both limiting energies are consequently
$\mathcal E(\K,\alpha\K)$. It follows that the original energy is
constant, and hence $v\equiv\K$. This proves the proposition.
\end{proof}

We can now complete the proof of the main theorem.

\begin{proof}[Proof of Theorem~\ref{thm:classification}]
Suppose that the origin is nonremovable. By
Proposition~\ref{prop:removability}, $u$ satisfies
\eqref{eq:two-sided-local}. Given an arbitrary sequence $r_j\downarrow0$,
set $u_j(y)=r_j^\alpha u(r_jy)$. The annular compactness in
Section~\ref{sec:prelim} gives, after passage to a subsequence,
\[
 u_j\longrightarrow U\quad\text{in }C^1_{\loc}(\R^n\setminus\{0\}),
\]
where $U>0$ solves $-\Delta_pU=U^q$ in the punctured whole space.
The two-sided bounds pass to the limit on every compact annulus,
so Proposition~\ref{prop:global-rigidity} gives
$U(y)=\K|y|^{-\alpha}$.

Every sequence of scales therefore has a subsequence with this same
limit. By compactness, the full family converges as in
\eqref{eq:C1limit}. Restricting to the unit sphere gives
\eqref{eq:precise}.
Since $\alpha>0$ and $\K>0$, this alternative is incompatible
with removability. The proof is complete.
\end{proof}

\noindent M.~Xu

\noindent School of Mathematical Sciences\\
Fudan University, Shanghai 200433, China \\
Email: \href{mailto:meiqing_xu@fudan.edu.cn}{\nolinkurl{meiqing_xu@fudan.edu.cn}}

\bigskip

\noindent H.~Yang

\noindent School of Mathematical Sciences\\
Shanghai Jiao Tong University, Shanghai 200240, China \\
Email: \href{mailto:hui-yang@sjtu.edu.cn}{\nolinkurl{hui-yang@sjtu.edu.cn}}

\end{document}